\newif\ifWorkingVersion
\WorkingVersionfalse  

\documentclass[12pt]{article}
\usepackage{amssymb} 

\newcommand{\beq}{\begin{eqnarray}}
\newcommand{\eeq}{\end{eqnarray}}
\newcommand{\beqn}{\begin{eqnarray*}}
\newcommand{\eeqn}{\end{eqnarray*}}

\newcommand{\be}{\begin{equation}}
\newcommand{\ee}{\end{equation}}
\newcommand{\bea}{\begin{eqnarray}}
\newcommand{\eea}{\end{eqnarray}}
\newcommand{\bean}{\begin{eqnarray*}}
\newcommand{\eean}{\end{eqnarray*}}
\newtheorem{defn}{Definition}
\newtheorem{thm}{Theorem}

\newtheorem{lemma}{Lemma}
\def\proof{\noindent{\bf Proof:}\quad}

\def\eps{\epsilon}
\def\Domain{{\{z:|z|\le \r+\e\}-[\r,\r+\e]}}

\def\a{{\alpha}}
\def\l{\lambda}

\def\si{\sigma}

\def\b{{\beta}}

\def\r{{\rho}}

\def\xh{{\hat x}}
\def\yb{{\hat y}}
\def\yh{{\hat y}}

\def\e{\epsilon}
\def\Qb{{\hat Q}}
\def\Qh{{\hat Q}}
\def\QS{Q^{\star}}
\def\xS{x^{\star}}
\def\yS{y^{\star}}
\def\ys{y^{\star}}

\def\pd{{\partial}}
\def\Mh{{\hat M}}

\def\O{{\tilde O}}

\def\fw{{\rm fw}}
\def\ew{{\rm ew}}
\def\Do{{D}}
\def\tO{{\hat t}}
\def\BOp{\phantom{x}}  
\def\net#1{{\cal #1}}  
\def\graph{\net}
\def\map{\graph}
\def\cf#1#2{\null_{#2}#1}  
\def\cf#1#2{#1^{[#2]}\null}  

\begin{document}

\ifWorkingVersion
\overfullrule=2pt  
\pagestyle{myheadings} 
\markboth={{\rm\jobname.tex\hfill \scriptsize(\LaTeX ed on \today)\hfill}}  
\fi
\def\Note#1{\ifWorkingVersion{\marginpar{\raggedright\tiny{#1}}}\fi\ignorespaces}

\title {Asymptotic enumeration of labelled graphs by genus}
\author{Edward A. Bender\\
Department of Mathematics\\
University of California at San Diego\\
La Jolla, CA 92093-0112, USA\\
\and
Zhicheng Gao\thanks{Research supported by NSERC}\\
 School of Mathematics and Statistics\\ Carleton
University\\ Ottawa Canada K1S 5B6 }
\date{}
\maketitle

\begin{abstract}
We obtain asymptotic formulas for the number of rooted
2-connected and 3-connected surface maps on an orientable surface
of genus $g$ with respect to vertices and edges simultaneously.
We also derive the bivariate version of the large face-width
result for random 3-connected maps.
These results are then used to derive asymptotic formulas
for the number of labelled $k$-connected graphs of orientable
genus $g$ for $k\le3$.
\end{abstract}

\section{Introduction}

The exact enumeration of various types of maps on the sphere (or,
equivalently, the plane) was carried out by Tutte~\cite{Tut62,
Tut62a, Tut63} in the 1960s via his device of rooting. (Terms in
this paragraph are defined below.) Building on this, explicit
results were obtained for some maps on low genus surfaces, e.g., as
done by Arqu\'es on the torus~\cite{Arq87}. Beginning in the 1980s,
Tutte's approach was used for the asymptotic enumeration of maps on
general surfaces~\cite{BC86, BW88, BCR93}. A matrix integral
approach was initiated by $'$t Hooft (see \cite{LZ04}). The
enumerative study of graphs embeddable in surfaces began much more
recently. Asymptotic results on the sphere were obtained
in~\cite{BGW02, McD08, GN09} and cruder asymptotics for general
surfaces in~\cite{McD08}. In this paper, we will derive asymptotic
formulas for the number of labelled graphs on an orientable surface
of genus $g$ for the following families: 3-connected and 2-connected
with respect to vertices and edges, and 1-connected and all with
respect to vertices. Along the way we also derive results for
2-connected and 3-connected maps with respect to vertices and edges.
The result for all graphs as well as various parameters for these
graphs was announced earlier by Noy~\cite{NoyAofA} and uploaded
in~\cite{CFGMN10}.

\begin{defn}[Maps and Embeddable Graphs]

A {\em map} $\map M$ is a connected graph $\graph G$ embedded in a surface $\Sigma$
(a closed 2-manifold) such that all components of $\Sigma-\graph G$
are simply connected regions, which are called {\em faces}.
$\graph G$ is called the underlying graph of $\map M$, and is denoted by $G(\map M)$.
Loops and multiple edges are allowed in $\graph G$.
\begin{itemize}
\item
A map is {\em rooted} if an edge is distinguished together with
a direction on the edge and a side of the
edge.
$$\mbox{In this paper, all maps are rooted and unlabeled.}$$
\item
A graph without loops or multiple edges is {\em simple}.
\item
A graph $\graph G$ is {\em embeddable in a surface} if it can be drawn on
the surface without edges crossing.
\item
A graph has (orientable) genus $g$ if it is embeddable in an
orientable surface of genus $g$ and none of smaller genus.
\end{itemize}
\end{defn}

\begin{defn}[Generating Functions for Maps and Graphs]~
Let $\Mh_g(n,m;k)$ be the number of (rooted, unlabeled)
$k$-connected maps with $n$ vertices and $m$ edges, on an
orientable surface of genus $g$.
Let $G_g(n,m;k)$ be the number of (vertex) labelled, simple,
$k$-connected graphs with $n$ vertices
and $m$ edges, which are embeddable in an orientable surface of
genus $g$. Let $G_g(n;k)=\sum_{m}G_g(n,m;k)$, the number of
labelled, simple,\break
$k$-connected graphs with $n$ vertices.
Let
$$
\Mh_{g,k}(x,y) = \sum_{n,m} \Mh_g(n,m;k)x^ny^m ~~\mbox{and}~~
G_{g,k}(x,y) = \sum_{n,m} G_g(n,m;k)(x^n/n!)y^m.
$$
\end{defn}

In the following theorem, $\rho(r)$ and $A_g(r)$ have the same
definition in terms of $r$, but the definition of $r$ varies.

\begin{thm}[Maps on Surfaces] \label{Thm:maps}
Define
\beqn
\r(r)&=&\frac{r^3(2+r)}{1+2r},\\
A_g(r)&=&\frac{1}{2\sqrt{\pi}}\frac{r^6(2+r)^{3/2}}{(1+2r)^2}
\left(\frac{12(1+r)^3(1+2r)^4}{r^{12}(2+r)^5}\right)^{g/2}t_g,
\eeqn
where $t_g$ is the map asymptotics constant defined in \cite{BC86}.
For $k=1,2,3$, there are algebraic functions $r=r_k(m/n)$, $C_k(r)$,
and $\eta_k(r)$ such that for any fixed $\e>0$ and fixed genus $g$
$$
\Mh_g(n,m;k) ~\sim~ C_k(r)A_g(r)(2+r)^{(k-1)(5g-3)/2}\,n^{5g/2-3}
\r(r)^{-n}\eta_k(r)^{-m},
$$
uniformly as $n,m\to \infty$
such that $r_k(m/n)\in [\e,1/\e]$. The relevant functions are as
follows:
\begin{itemize}
\item[(i)]
$\displaystyle r=r_1(m/n) \mbox{ satisfies }
\frac{(1+r)(1+r+r^2)}{r^2(2+r)}=\frac{m}{n},~~
\eta_1(r)=\frac{1+2r}{4(1+r+r^2)^2}$~~and
$$
C_1(r)~=~(2+r)\sqrt{\frac{1+r+r^2}{(1+2r)(4+7r+4r^2)}};
$$
\item[(ii)]
$\displaystyle r_2(m/n)=\frac{1}{m/n-1},~~
\eta_2(r)=\frac{4}{(1+2r)(2+r)^2_{\vphantom{\Bigm|}}} ~~\hbox{and}~~
C_2(r) = \frac{1}{\sqrt{(1+2r)(2+r)}}$;
\item[(iii)]
$\displaystyle r_3(m/n)=\frac{3-m/n}{2(m/n)-3},~~
\eta_3(r)=\frac{3}{4r(2+r)_{\vphantom{\Bigm|}}}, ~~\hbox{and}~~
C_3(r)=\frac{1}{\sqrt{r(2+r)^3}}$.
\end{itemize}
\end{thm}

\begin{thm}[Embeddable Graphs] \label{Thm:graphs}
For the ranges of $m$ and $n$ considered here,
the number of graphs embeddable in an orientable
surface of genus $g$ is asymptotic to the number
of such graphs of orientable genus $g$.
\begin{itemize}
\item[(i)]
{\bf(3-connected)}
For any fixed $\e>0$ and genus $g$,
$$
\frac{G_g(n,m;3)}{n!} ~\sim~ \frac{\Mh_g(n,m;3)}{4m}
$$
uniformly as $n,m\to \infty$ such that $\frac{m}{n}\in
[(3/2)+\e,3-\e]$.
\item[(ii)]
{\bf(2-connected)} Let $\a(t),\b(t)$, $\rho_2(t)$, $\l_2(t)$,
$\mu(t)$ and $\si(t)$ be functions of $t$ defined in
Section~\ref{Sec:3to2} (see also \cite{BGW02}).
Let
$$
B_g(t)=\left(\frac{8}{9(1+t)(1-t)^6}\left(\frac{\b(t)}{\a(t)}\right)^{5/2}\right)^{g-1}.
$$
Fix $\e>0$ and genus $g$. Let $0<t<1$ satisfy $\mu(t)=m/n$.
Then
$$
\frac{G_g(n,m;2)}{n!} ~\sim~ \frac{B_g(t)t_g}{4\si(t)\sqrt{2\pi}}
n^{5g/2-4}\rho_2(t)^{-n}\l_2(t)^{-m}
$$
 uniformly as $n,m\to \infty$ such that $m/n \in
[1+\e,3-\e]$.
\item[(iii)]
{\bf(vertices only)} For $0\le k\le 3$ and fixed $g$, there are
positive constants $x_k$, $\a_k$ and $\b_k$ such that
$$
\frac{G_g(n;k)}{n!} ~\sim~ \a_k\b_k^g\,t_g\,n^{5g/2-7/2}x_k^{-n},
$$
where
$$\matrix{
x_3\doteq 0.04751,\hfill & x_2\doteq 0.03819,\hfill & x_1\doteq
0.03673,\hfill & x_0=x_1,\hfill\cr \b_3\doteq 1.48590\cdot
10^{5},\hfill & \b_2\doteq 7.61501\cdot 10^4. \hfill & \b_1\doteq
6.87242\cdot 10^4,\hfill & \b_0=\b_1,\hfill\cr
 \a_3=\frac{1}{4\b_3},\hfill & \a_2=\frac{1}{4\b_2},\hfill &
 \a_1=\frac{1}{4\b_1},\hfill & \a_0\doteq 3.77651\cdot 10^{-6}.\hfill\cr }
$$
More accurate values of these constants can be computed by
using the formulas in those sections where the theorem
is proved.
\end{itemize}
\end{thm}

\smallskip\noindent
{\bf Remark ($t_g$)}.
It is known \cite{GLM08} that
$$
t_g ~=~ \frac{-a_g}{2^{g-2}\Gamma\left(\frac{5g-1}{2}\right)}
$$
where $a_0=1$ and, for $g>0$,
\be\label{ag}
a_g ~=~ \frac{(5g-4)(5g-6)}{48}a_{g-1}-\frac{1}{2}\sum_{h=1}^{g-1}a_{h}a_{g-h}.
\ee
Hence all the numbers in Theorems~\ref{Thm:maps} and~\ref{Thm:graphs}
can be computed efficiently to any desired accuracy for any given $g$ and $r$.

\smallskip\noindent
{\bf Remark (Sharp Concentration)}.
As noted in Comment~4 of Section~\ref{Sec:Tools}, our methods for obtaining
bivariate results show that the number of edges is sharply concentrated.
To find the mean number of edges asymptotically,
set $\eta_k(r)=1$ in Theorem~\ref{Thm:maps},
$\eta_3(r)=1$ in Theorem~\ref{Thm:graphs}(i),
and $\lambda_2(t)=1$ in Theorem~\ref{Thm:graphs}(ii).
For $r$ the asymptotic value of the mean is then the value of $m$ for
which $r(m/n)$ has that value of $r$; for $t$ it is simply $\mu(t)n$.

\medskip
The paper proceeds as follows.

\begin{description}
\item[Section~\ref{Sec:Quad}]
Maps on a fixed surface were enumerated in
\cite{BCR93} with respect to vertices and faces.
We convert this result to quadrangulations and then obtain results for
other types of quadrangulations.

\item[Section~\ref{Sec:Tools}:]
We recall a local limit theorem and discuss some analytic methods used in
subsequent sections.

\item[Section~\ref{Sec:Maps}:]
We then apply the techniques in \cite{BW88} and \cite{BGRW96} to obtain
asymptotics for generating functions for $k$-connected maps, proving
Theorem~\ref{Thm:maps}.
The calculations for $A_g(r)$ are postponed to Section~\ref{Sec:Ag(r)}.

\item[Section~\ref{Sec:FW}:]
Applying the techniques in \cite{BGR94}, we show that almost all
3-connected maps have large face-width when counted by vertices and
edges. Hence almost all 3-connected graphs of genus $g$ have a
unique embedding~\cite{RV90}. This leads to Theorem~\ref{Thm:graphs}
for 3-connected graphs.

\item[Section~\ref{Sec:3to2}:]
Using the construction of 2-connected graphs from 3-connected
graphs and polygons as in~\cite{BGW02} we obtain Theorem~\ref{Thm:graphs}
for 2-connected graphs.

\item[Sections~\ref{Sec:2to1} and \ref{Sec:1toAll}:]
We obtain Theorem~\ref{Thm:graphs} for 1-connected graphs from
the 2-connected result and for all graphs from 1-connected by methods like
those in~\cite{GN09}.

\item[Section~9:] We derive the expression for $A_g(r)$ in terms of
$t_g$.

\item[Section~10:] We make some comments on the number of labeled graphs of a given
nonorientable genus.
\end{description}

\section{Enumerating Quadrangulations}\label{Sec:Quad}

We begin with some definitions:
\begin{defn}[Cycles]
A {\em cycle} in a map is a simple closed curve consisting of edges of the map.
\begin{itemize}
\item A cycle is called a {\em $k$-cycle} if it contains $k$ edges.
\item
A cycle is called {\em separating} if deleting it separates the underlying graph.
\item
A cycle is called {\em facial} if it bounds a face of the map.
\item
A cycle is called {\em contractible} if it is homotopic to a point,
otherwise it is called {\em non-contractible}.
\item
A contractible cycle in a nonplanar map separates the map into a
planar piece and a nonplanar piece.
The planar piece is called the {\em interior} of the cycle and we
also say that the cycle {\em contains} anything that is in its interior.
Since we usually draw a planar map such that the root face is the
unbounded face, we define the interior of a cycle in a planar map to
be the piece which does not contain the root face.
\item
A 2-cycle or 4-cycle is called {\em maximal\/} ({\em minimal\/})
if it is contractible and its interior is maximal (minimal).
\end{itemize}
\end{defn}
\begin{defn}[Widths]
The {\em edge-width} of a map $\map M$, written $\ew(\map M)$, is the
length of a shortest non-contractible cycle of $\map M$.
The {\em face-width} (also called {\em representativity} of $\map M$,
written $\fw(\map M)$, is the minimum of $|G(\map M)\cap C|$ taken over all
non-contractible closed curves $C$ on the surface.
\end{defn}
\begin{defn}[Quadrangulations]
A quadrangulation is a map all of whose faces have degree 4.
\begin{itemize}
\item
A {\em bipartite quadrangulation} is a quadrangulation whose
underlying graph is bipartite.
(All quadrangulations on the sphere are bipartite, but
those on other surfaces need not be.)
\item
A quadrangulation is {\em near-simple} if it has no
contractible 2-cycles and no contractible nonfacial 4-cycles.
\item
A quadrangulation is {\em simple} if it has no 2-cycles and all 4-cycles are
facial.
\end{itemize}
\end{defn}

\noindent The following lemma, contained in \cite{BW88} and \cite{BGRW96},
connects maps with bipartite quadrangulations.
\begin{lemma} \label{Lemma:bijection}
By convention, we bicolor a bipartite quadrangulation so that the head of the
root edge is black.
 There is a bijection $\phi$ between rooted maps and
rooted bipartite quadrangulations, such that the following hold.
\begin{itemize}
\item[(a)]
$\fw(\map M)=\ew(\phi(\map M))/2$.
\item[(b)]
$\map M$ has $n$ vertices and $m$ edges if and only if $\phi(\map M)$
has $n$ black vertices and $m$ faces.
\item[(c)]
$\phi(\map M)$ has no 2-cycle  implies $\map M$ is 2-connected which implies
$\phi(\map M)$ has no contractible 2-cycle.
\item[(d)]
$\phi(\map M)$ is simple implies $\map M$ is 3-connected which implies
$\phi(\map M)$ is near-simple.
\end{itemize}
\end{lemma}

In this section we enumerate quadrangulations with no contractible
2-cycles and near-simple quadrangulations.
Except that black vertices were not counted, this is done in \cite{BGRW96}.
In what follows, we reproduce that argument nearly verbatim, adding
a second variable to count black vertices.

We define the generating functions $Q_g(x,y)$, $\Qb_g(x,y)$ and
$\QS_g(x,y)$ as follows.
$$
Q_g(x,y)=\sum_{i,j\ge1}Q(i,j;g)x^{i-1}y^{j}
$$
where $Q(i,j;g)$ is the number of (rooted, bicolored)
quadrangulations with $i$ black vertices and $j$ faces on an orientable
surface of genus $g$.
Similarly define $\Qb_g(x,y)$ for quadrangulations without contractible
2-cycles and $\QS_g(x,y)$ for near-simple quadrangulations.

By Lemma \ref{Lemma:bijection}, we have
\be \label{QMrel}
Q_g(x,y)=x^{-1}\Mh_{g,1}(x,y)-\delta_{0,g},
\ee
where the Kronecker delta occurs because of the convention that
counts a single vertex as a map on the sphere.

In \cite{BCR93} the generating function $\Mh_g(u,v)$ counts maps by
vertices and faces. Thus \be\label{Mhats} \Mh_{g,1}(x,y) ~=~
y^{2g-2}\Mh_g(xy,y). \ee It is known \cite{Arq87,BCR93} that
$\Mh_0(xy,y)=\frac{rs}{(1+r+s)^3}$  where $r(x,y)$ and $s(x,y)$ are
power series uniquely determined by \be \label{xyrs} x ~=~
\frac{r(2+r)}{s(2+s)}~~\mbox{and}~~y ~=~ \frac{s(2+s)}{4(1+r+s)^2}.
\ee Thus \be\label{Q0} Q_0(x,y) ~=~ \frac{4(1+r+s)}{(2+r)(2+s)}-1
~=~ \frac{2r+2s-rs}{(2+r)(2+s)}, \ee and \be \label{Jacob} \matrix{
\displaystyle \frac{\pd r}{\pd x} ~=~
\frac{s(2+s)(1+r+rs)}{2(1-rs)},\hfill&~~ \displaystyle \frac{\pd
r}{\pd y} ~=~ \frac{2r(2+r)(1+s)(1+r+s)^3}{s(2+s)(1-rs)},\cr
\displaystyle \frac{\pd s}{\pd x} ~=~
\frac{s^2(2+s)^2}{2(1-rs)},\hfill&~~ \displaystyle \frac{\pd s}{\pd
y} ~=~ \frac{2(1+r)(1+r+s)^3}{1-rs}.\hfill} \ee

Throughout the rest of the
paper, we use $N(\e)$ to denote the set
$$
N(\e) ~=~ \{re^{i\theta}:\e \!\le\! r \!\le\! 1/\e,\; |\theta|\le \e\}.
$$

\begin{thm}[Quadrangulations]\label{Thm:Quad}
Fix $g>0$ and let $q(x,y)$ be any of $Q_g(x,y)$, $\Qh_g(x,y)$ and $\QS_g(x,y)$.
The values of $x$ and $y$ are parameterized by $r$ and $s$ in the following
manner.
\begin{itemize}
\item[(i)]
For all (bipartite) quadrangulations ($q=Q_g$), $x$ and $y$ are given by (\ref{xyrs}).
\item[(ii)]
For no contractible 2-cycles ($q=\Qh_g$), $x$ is given by (\ref{xyrs}) and
$\displaystyle y = \frac{4s}{(2+s)(2+r)^2}$.
\item[(iii)]
For near simple ($q=\QS_g$), $x$ is given by (\ref{xyrs}) and
$\displaystyle y = \frac{s(4-rs)}{4(2+r)}$.
\end{itemize}
The following are true.
\begin{itemize}
\item[(a)]
The function $q(x,y)$ is a rational function of $r$ and $s$
and hence an algebraic function of $x$ and $y$.
\item[(b)]
If $r$ and $s$ are positive reals such that $rs=1$, then $(x,y)$ is in the
singular set of $q(x,y)$.
\item[(c)]
If $(x',y')$ is another singularity of $q$, then either $|x'|>x$ or $|y'|>y$.
\item[(d)]
Let $\rho(r)=\frac{r^3(2+r)}{1+2r}$, the value of $x$ on the
singular curve $rs=1$, and let $y$ be its value on the singular
curve at $r$. Fix $\eps>0$ and $g>0$. Uniformly for $r\in N(\eps)$
\be\label{Qasymp} xq(x,y) ~\sim~
C(r)\left(1-\frac{x}{\rho(r)}\right)^{(3-5g)/2} \ee as
$x\to\rho(r)$,
$$
C(r) = \cases{ \displaystyle\sqrt{\frac{\pi}{3(1+r)}}\,
\frac{(1+r+r^2)A_g(r)\,\Gamma\!\left(\frac{5g-3}{2}\right)}{r^2}
_{\vphantom{\Bigr|}} &for $q=Q_g$,\cr
\displaystyle\sqrt{\frac{\pi}{3(1+r)}}\,
\frac{A_g(r)\,\Gamma\!\left(\frac{5g-3}{2}\right)}{r}(2+r)^{(5g-3)/2}
\vphantom{\Biggr|^{\Bigr|}_{\Bigr|}} &for $q=\Qh_g$,\cr
\displaystyle\sqrt{\frac{3\pi}{1+r}}\,
\frac{A_g(r)\,\Gamma\!\left(\frac{5g-3}{2}\right)}{(2+r)(1+2r)}(2+r)^{5g-3}
&for $q=\QS_g$,}
$$
and some function $A_g(r)$ whose
value is determined in Section~\ref{Sec:Ag(r)}.
\end{itemize}
\end{thm}

\proof
Theorem~3 of \cite{BCR93} shows that $\Mh_g(x,y)$ of that paper is a
rational function of $r$ and $s$ and hence algebraic when $g>0$.
(The theorem contains the misprint $9>0$ which should be $g>0$.)
Use~(\ref{QMrel})--(\ref{Q0}) to establish~(a) for $Q_g$.

We now derive equations for $\Qh$ and $\QS$ based on $Q$.
This will easily imply (a) for $\Qh$ and $\QS$.

It is important to note that, in any quadrangulation, all
maximal 2-cycles have disjoint interiors, and that, in any nonplanar
quadrangulation without contractible 2-cycles, all maximal 4-cycles
have disjoint interiors. (This is simpler than the planar case
\cite[p.\thinspace260]{MS68}.) Therefore, we can close all maximal 2-cycles
in quadrangulations to obtain quadrangulations without contractible
2-cycles and remove the interior of each maximal contractible
4-cycle to obtain near-simple quadrangulations. The process can be
reversed and used to construct quadrangulations from near-simple
quadrangulations.

\medskip\noindent
{\bf Enumerating $\Qb_g(x,y)$:} The following argument is
essentially from \cite{BGRW96}, by paying extra attention to the
number of black vertices. All quadrangulations of genus $g>0$ can be
divided into two classes according as the root face lies in the
interior of some contractible 2-cycle or not.

For any quadrangulation in the first class, let $C$ be the minimal
contractible 2-cycle containing the root face.
Cutting along $C$, filling
holes with disks and closing those two 2-cycles, we obtain a general
quadrangulation of genus $g$ and a planar quadrangulation with a
distinguished edge. Taking the latter quadrangulation and cutting
along all its maximal 2-cycles and closing as before gives a
quadrangulation without contractible 2-cycles, together with a set
of planar quadrangulations extracted from within the maximal
2-cycles.
Remembering that $y$ counts faces and that the number of edges is
twice the number of faces, it follows that the generating function
for the first class is
$$
\frac{Q_g(x,y)}{1+Q_0(x,y)}\,
\frac{2\yb \;\partial\,\Qb_0(x,\yb)}{\partial \yb},
$$
where \be \label{yb} \yb ~=~ y(1+Q_0(x,y))^2 ~=~
\frac{4s}{(2+s)(2+r)^2}. \ee For any quadrangulation in the second
class, closing all maximal contractible 2-cycles gives
quadrangulations without contractible 2-cycles. Thus the generating
function for this class is $\Qb_g(x,\yb)$. For the planar case, only
the second class applies and so
 \be \label{Qb0}
\Qb_0(x,\yb)=Q_0(x,y). \ee Combining the two classes when $g>0$, we
have
$$
Q_g(x,y) ~=~ \Qb_g(x,\yb)+\frac{Q_g(x,y)}{1+Q_0(x,y)}\,
\frac{2\yb\; \partial\,\Qb_0(x,\yb) }{\partial \yb}.
$$
It follows that \beq \label{Qb} \Qb_g(x,\yb) ~=~
\left(1-\frac{2\yb}{1+Q_0(x,y)}\,
\frac{\partial\,\Qb_0(x,\yb)}{\partial\yb}\right)Q_g(x,y) \eeq for
$g>0$. Note that \beq\label{Qb factor} 1-\frac{2\yb}{1+Q_0(x,y)}\,
\frac{\partial\,\Qb_0(x,\yb)}{\partial\yb} ~=~ \frac{1}{1+r+s} \eeq
and so is bounded on the singular curve when $r$ is near the
positive real axis.

\noindent {\bf Enumerating $\QS_g(x,y)$:}
We now use a similar argument to derive $\QS_g(x,\ys)$ from
$\Qb_g(x,\yb)$ when $g>0$.
For any quadrangulation without contractible 2-cycles, let $C$ be the
maximal contractible 4-cycle containing the root face. Cutting along
$C$ and filling holes with disks,  we obtain
\begin{enumerate}
\item a planar quadrangulation which has no 2-cycles and
has a distinguished face other than the root face, and
\item a quadrangulation of genus $g$ which, after the
removal of the interiors of all maximal 4-cycles, gives a
near-simple quadrangulation.
\end{enumerate}
Note that
\be \label{ys}
\ys ~=~ \frac{\Qb_0(x,\yb)-x\yb-\yb}{x\yb} ~=~ \frac{s(4-rs)}{4(2+r)}
\ee
enumerates planar quadrangulations having at least one interior face
and having no 2-cycles such that $x$ marks the number of black
vertices minus 2 and $\yb$ marks the number of non-root faces. It
follows from the construction that
$$
\frac{\Qb_g(x,\yb)}{\yb} ~=~ \frac{\QS_g(x,\ys)}{\ys}
\frac{\pd \ys}{\pd\yb}.
$$
which gives \beq \label{QS} \QS_g(x,\ys) ~=~ \frac{\ys}{\pd
\ys/\pd\yb}\frac{\Qb_g(x,\yb)}{\yb} ~=~
\frac{4-rs}{(2+s)(2+r)(1+r+s)}Q_g(x,y). \eeq
This completes the proof of Theorem~\ref{Thm:Quad}(a).

\medskip\noindent {\bf Singularities:}
These must arise from poles due to the vanishing of the denominator
of $q(x,y)$ or from branch points caused by problems with the Jacobian
$\frac{\partial(x,y)}{\partial(r,s)}$.
For the former, it can be seen from (\ref{Qb}) and (\ref{QS}) that either
$1+r+s=0$ or $2+r=0$ or $2+s=0$.
By (\ref{xyrs}), each of these implies that either $x$ or $y$ vanishes
or is infinite, which do not matter since the radius of convergence
is nonzero and finite.
Using the formulas in Theorem~\ref{Thm:Quad}, one can
compute Jacobians.
One finds that the only singularity that matters is $1-rs=0$.

Conclusion (c) follows for $Q$ from \cite{BCR93}.
We now consider $\Qh$ and $\QS$.
Suppose
\begin{itemize}
\item

$x$ and $y$ are positive reals on the singular curve,
\item
$x'$ and $y'$ are on the singular curve,
\item
$|x'|\le x$ and $|y'|\le y$.
\end{itemize}
To prove (c) it suffice to show that $x'=x$ and $y'=y$.
Since we are dealing with generating functions with nonnegative coefficients, no singularity
can be nearer the origin the that at the positive reals.
Hence $|x'|=x$ and $|y'|=y$.
As was done in \cite{BR84}, one easily verifies that on the singular curve $rs=1$ one has
\be\label{QS singular}
16x'y'^2\Bigl(16(y'+1)(x'y'+1)+2\Bigr) ~=~ 27
\ee
for $\QS$.
Taking absolute values in this equation one easily finds that $|y'+1|=|y+1|$
and $|x'y'+1|=|xy+1|$.
Thus $y'=y$ and $x'y'=xy$ and we are done.
For $\Qh$, a look at the equations for $x$ and $y$ on the singular curve shows that
we need only replace $y'$ in (\ref{QS singular}) with
$(3/4)(y'/4x')^{1/3}$ and argue as for $\QS$.
This completes the proof of~(c).

\medskip\noindent {\bf Asymptotics:}
We now turn to (d). The case $q=Q_g$ is contained implicitly
in~\cite{BCR93} for some function $A_g(r)$.

We now use (\ref{Qb}) to derive the singular expansion for
$\Qh_g(\xh,\yh)$ at $\xh=\r(r)$ where $r$ is determined by
$\yh=\eta_2(r)$. It is important to note that, with $\yh$ fixed,
(\ref{yb}) defines $y$ as an analytic function in $x=\xh$. Thus in
(\ref{Qasymp}), with $q(x,y)=Q_g(x,y)$, we should treat $r$ as a
function in $y$ and consequently as a function in $x$. Using
implicit differentiation, we obtain from (\ref{yb}) and
(\ref{Jacob}) that \be \label{dydx} \frac{dy}{dx}=-\frac{\pd \yh/
\pd x }{\pd \yh/\pd y}
 ~=~ -\frac{(\pd\yh/ \pd r)(\pd r/\pd x)+(\pd \yh/ \pd s)(\pd s/\pd x)}
{(\pd\yh/ \pd r)(\pd r/\pd y)+(\pd \yh/ \pd s)(\pd s/\pd y)}
~=~ \frac{-s^2(2+s)^2}{4(2+r)(1+r+s)^3}.
 \ee
Hence $$\frac{d }{dx}\left(1-\frac{x}{\r(r)}\right)
~=~ \frac{-1}{\r(r)}+\frac{x}{\r^2(r)}\frac{d\r}{dx}
~=~ \frac{-1}{\r(r)}+\frac{x}{\r^2(r)}\frac{\r'(r)}{\eta'_1(r)}\frac{dy}{dx}.
$$
Using (\ref{dydx}) and the expressions for $\r(r)$ and $\eta_1(r)$
given in Theorem~1, we obtain
$$
\left.\frac{d}{dx}\left(1-\frac{x}{\r(r)}\right)\right|_{x=\r(r),s=1/r}
\!=~ \frac{-1}{\r(r)(2+r)},
$$
and hence
$$1-\frac{x}{\r(r)}
~\sim~ \frac{-1}{\r(r)(2+r)}(\xh-\rho(r))
~=~ \frac{1}{2+r}\left(1-\frac{\xh}{\rho(r)}\right).
$$

 Substituting this into (\ref{Qasymp}), we obtain
$$
\left(1-\frac{x}{\r(r)}\right)^{(3-5g)/2} \!\sim~
(2+r)^{(5g-3)/2}\left(1-\frac{\xh}{\r(r)}\right)^{(3-5g)/2},
$$
as $\xh\to \rho(r)$ for each fixed $\yh$. The factor (\ref{Qb
factor}) can simply be evaluated at $s=1/r$ since it converges to a
constant. This establishes (\ref{Qasymp}) for $\Qh_g(\xh,\yh)$.

Expansion (\ref{Qasymp}) for $\QS_g(\xS,\yS)$ can be obtained
similarly using (\ref{QS}). We note that fixing $\yS$ defines $y$,
and hence $\r(r)$, as a function of $x=\xS$. Using (\ref{ys}) and
(\ref{Jacob}), we obtain
$$
1-\frac{x}{\r(r)} ~\sim~ \frac{1}{(2+r)^2}\left(1-\frac{\xS}{\rho(r)}
\right),
$$
as $\xS\to \rho(r)$ for each fixed $\yS$.

This completes the proof of the theorem, except for the formula for $A_g(r)$
which will be derived in Section~\ref{Sec:Ag(r)}.

\section{Some Technical Lemmas}\label{Sec:Tools}

The following lemma is the essential tool for our asymptotic estimates.
It is based on the case $d=1$ of~\cite[Theorem 2]{BR83},
 from which it follows immediately.
\begin{lemma}\label{Lemma:LLT}
Suppose that $a_{n,k}\ge0$. Define $a_n(v) =\sum_{k}a_{n,k}v^k$ and
$a(u,v) = \sum_n a_n(v)u^n$. Let $R(c)$ be the radius of convergence
of $a(u,c)$. Suppose that $I$ is a closed subinterval of
$(0,\infty)$ on which $0<R<\infty$. For $v\in I$ define
$$
\mu(v) = \frac{-d\,\log\rho(v)}{d\,\log v},~~
\sigma^2(v) = \frac{-d^2\,\log\rho(v)}{(d\,\log v)^2},~~
K_n=\{n\mu(v)\mid v\in I\}\cap{\mathbb Z}
$$
and $N(I,\delta)=\{z\mid |z|\in I~\mbox{and}~|\arg z|<\delta\}$.
Suppose there are $f(n)$, $g(v)$
and $\r(v)$ such that in $N(I,\delta)$
\begin{itemize}
\item[(a)] $a_n(v) \sim f(n)g(v)\r(v)^{-n}$ uniformly as $n\to\infty$;
\item[(b)] $g(v)$ is uniformly continuous;
\item[(c)] $\r(v)\ne 0$ has a uniformly continuous third derivative;
\item[(d)] $\sigma^2(v)>0$ for $v>0$.
\end{itemize}
Suppose also that
\begin{itemize}
\item[(e)] $R(c)>R(|c|)$ whenever $c\ne|c|\in I$.
\end{itemize}
Then, as $n\to\infty$, we have, uniformly for $k\in K_n$,
$$
a_{n,k} ~\sim~ \frac{a_n(v)v^{-k}}{\sqrt{2\pi n\sigma^2(v)}},
$$
where $v\in I$ is given by $k/n=\mu(v)$.
\end{lemma}
Of course $|\r(v)|$ is simply the radius of convergence $R(v)$ and
$\r(v)=R(v)$ when $v\in I$.

We now make some comments on applying this lemma.
We will generally use these ideas without explicit mention.

\smallskip\noindent
{\bf Comment 1}. There is the direct application. We can apply the
lemma to (\ref{Qasymp}) to obtain asymptotics. The only condition
that is not immediate is the verification that $\sigma^2(v)>0$ for
(d). This is a straightforward but somewhat tedious calculation.
Unless needed later, we omit the values of $\si^2(v)$ that we
compute.

\smallskip\noindent
{\bf Comment 2}.
There is the effect of adding and multiplying various $a(u,v)$,
all with the same $\r(v)$ (and hence $\mu(v)$) that satisfy the lemma.
The result will be a function that again satisfies the lemma with the same $\r(v)$.

To see this, note that the lemma is essentially a local limit theorem for
random variables where $\Pr(X_n\!=\!k)=a_{n,k}v^k/a_n(v)$ and
use~\cite[Lemma~5]{BRW83}.
We also need the observation that multiplying $a(u,v)$ by functions
with nonnegative coefficients and larger radii of convergence
results in a function having the same $\r(v)$ and so the lemma applies.
In fact, it suffices to simply evaluate the new function at the singularity
and multiply the resulting constant by $a(u,v)$.

\smallskip\noindent
{\bf Comment 3}. Condition (a) will follow if $a(u,v)$ is algebraic
and $a(u,s)$ has no other singularities on its circle of convergence
when $s\in I$. In general, condition (a) is established using the
``transfer theorem''~\cite[Sec.\thinspace VI.3]{FSbook}. Thus, for
example, Theorem~\ref{Thm:Quad}(a,c) implies
Lemma~\ref{Lemma:LLT}(a,e).

\smallskip\noindent
{\bf Comment 4}.
The values $n\mu(v)$ and $n\si^2(v)$ are asymptotic to the mean and variance of a random
variable $X_n(v)$ with $\Pr(X_n(v)\!=\!k) = a_{n,k}v^k/a_n(v)$.
Chebyshev's inequality then gives a sharp concentration result for $X_n(v)$ about its mean.
When this is applied to maps or graphs with $v=1$, it gives a sharp concentration for the
edges about the mean.
(The lemma is based on a local limit theorem, which could be used to give a sharper result.)

\medskip

Since we will be bounding coefficients of generating functions,
the following definition and lemma will be useful.
\begin{defn}[$\O$]
Let $A(x,y)$ and $B(x,y)$ be generating functions and let $B(x,y)$
have nonnegative coefficients. We write $A(x,y)=\O(B(x,y))$ if
there is a constant $K$ such that
$$
\left|[x^iy^j]\,A(x,y)\right|\le K[x^{i}y^j]\,B(x,y)~\hbox{ for all } i,j.
$$
\end{defn}

\begin{lemma} \label{Lemma:O(GF)}
Let $A(x,y)$, $B(x,y)$, $C(x,y)$, $D(x,y)$ and $H(x,y)$ be generating
functions, and $C(x,y)$, $D(x,y)$ and $H(x,y)$  have nonnegative
coefficients.
If $A(x,y)=\O(C(x,y))$ and $B(x,y)=\O(D(x,y))$, then
\begin{itemize}
\item[(i)] {\bf differentiation:}
$A_x(x,y)=\O(C_x(x,y))$ and $A_y(x,y)=\O(C_y(x,y))$;
\item[(ii)] {\bf integration:}
$\displaystyle \int_0^x A(x,y)dx=\O\left(\int_0^x C(x,y)dx\right)$ and
\hfill\break\phantom{{\bf integration:}}
$\displaystyle \int_0^y A(x,y)dy=\O\left(\int_0^y C(x,y)dy\right)$;
\item[(iii)] {\bf product:}
$A(x,y)B(x,y)=\O(C(x,y)D(x,y))$;
\item[(iv)] {\bf substitution:}
$A(H(x,y),y)=\O(C(H(x,y),y)$ and
\hfill\break\phantom{{\bf substitution:}}
$A(x,H(x,y))=\O(C(x,H(x,y))$\hfill\break
provided that the compositions as formal power series are well
defined.
\end{itemize}
\end{lemma}
The proof follows immediately from the definition of $\O$.

Obviously the definition of $\O$ and Lemma~\ref{Lemma:O(GF)}
can be stated for any number of variables.

We want to apply Lemma~\ref{Lemma:LLT} to $a(u,v)=A(u,v)+E(u,v)$ or
$a(u,v)=A(u,v)-E(u,v)$ when $A$ is a function and we know $E$ only
approximately.
Of course, this cannot be done directly since derivatives are involved.

The lemma will apply to $A(u,v)$ for $v\in
I$. We could attempt to estimate coefficients of $E(u,v)$ by some
crude method, but this fails because the order of growth of $E(u,v)$
is not sufficiently smaller than that of $A(u,v)$. What we will have
is that $E(u,v)=\O(F(u,v))$ where $F$ is a function built
from functions to which the lemma applies and which have dominant
singularities only where $A$ has them. Thus both functions have the
same $\r(v)$. Furthermore, the function $f(n)$ for $A$ grows
faster than the $f(n)$ for $F$. This is enough to show that
the coefficients of $F$ are negligible compared to those of $A$
because of Comment~2 above.
We will use these ideas without explicit mention when
considering error bounds.

\section{Proof of Theorem~\ref{Thm:maps}}\label{Sec:Maps}

The value of $A_g(r)$ in this section is simply the value assumed
in the proof of Theorem~\ref{Thm:Quad} in Section~\ref{Sec:Quad}.
The formula for $A_g(r)$ will be derived in Section~\ref{Sec:Ag(r)}.

For $g=0$ we find it easier to verify that the formulas in
Theorem~\ref{Thm:maps} agree with known results.
The $g=0$ case for general maps will follow when we use~\cite{BCR93} to
evaluate $A_g(r)$ in Section~\ref{Sec:Ag(r)}.
For maps with $i+1$ vertices and $j+1$ faces
the number of 2-connected planar maps equals~\cite{BT64}
$$
\frac{(2i+j-2)!\,(2j+i-2)!}{i!\,j!\,(2i-1)!\,(2j-1)!}
$$
and the number 3-connected planar maps is asymptotic to
$$
\frac{1}{3^5ij}{2i\choose j+3}{2j\choose i+3}
$$
uniformly as $\max(i,j)\to\infty$ \cite{BW88a}.
The verification of Theorem~\ref{Thm:maps} now requires only
some straightforward estimates of factorials and the fact that
$t_0=\frac{2}{\sqrt\pi}$.

We now assume $g>0$.

We derive the 1-connected case from Theorem~\ref{Thm:Quad}.
Lemma~\ref{Lemma:bijection} tells us that
$xQ_g(x,y)$ counts 1-connected maps by vertices and edges. Now apply
Theorem~\ref{Thm:Quad} and Lemma~\ref{Lemma:LLT}.
With $A_g(r)$ given by Theorem~\ref{Thm:Quad}, it follows that
$$
\Mh_g(n,m;1) ~\sim~
\frac{A_g(r)}{\si_1(r)\sqrt{2\pi}}n^{5g/2-3}\r(r)^{-n}\eta_1(r)^{-m}
$$
where \beqn \frac{m}{n}&=&\frac{-d\,\log\rho(r)}{d\,\log
\eta_1(r)} ~=~ \frac{(1+r)(1+r+r^2)}{r^2(2+r)}, \\
\si_1^2(r)&=&\frac{-d^2\,\log\rho(r)}{(d\,\log\eta_1(r))^2}
~=~ \frac{(4+7r+4r^2)(1+2r)(1+r+r^2)}{6r^4(2+r)^2(1+r)}. \eeqn This
gives Theorem~1(i). (Of course, we could also have
cited~\cite{BCR93}, but we need the derivation from
Theorem~\ref{Thm:Quad} so that we can evaluate $A_g(r)$ later.)

\medskip
Our proof for 2- and 3-connected maps uses
Lemma~\ref{Lemma:bijection} in connection with
Theorem~\ref{Thm:Quad} and Lemma~\ref{Lemma:LLT}.
We obtain upper and lower bounds from Lemma~\ref{Lemma:bijection}(c,d).
We show that Lemma~\ref{Lemma:LLT} can be applied to both bounds and that
the asymptotics are the same.

Upper bounds are provided by $\Qh$ and $\QS$.
These can be treated in the same manner as Theorem~\ref{Thm:maps}(i) was
derived from $Q$.
Let $E(x,y)$ be the errors in these upper bounds.
We handle $E(x,y)$ as discussed at the end of Section~\ref{Sec:Tools},
namely $E(x,y)=\O(F(x,y))$ where $F$ is well-behaved.
We now turn to $F(x,y)$.

\medskip
\noindent{\bf 2-Connected maps:}
We bound the quadrangulations counted by $\Qb_g(x,y)$ that have
non-contractible 2-cycles.
The argument is essentially the same as that
used in~\cite{BGRW96}. The only difference is that we keep track of
both the number of faces and the number of black vertices.

We first consider quadrangulations counted by $\Qb_g(x,y)$ which
contain a separating non-contractible cycle $C$ of length $2$.

Cutting through $C$ gives two near-quadrangulations. After closing
the resulting two 2-cycles, we obtain a rooted quadrangulation $\map Q_1$
with a distinguished edge, which has genus $0<j<g$, and another
rooted quadrangulation $\map Q_2$ with genus $g-j$. The quadrangulation
$\map Q_1$ may contain contractible 2-cycles which contain the
distinguished edge $d$ in its interior.
Hence $\map Q_1$ is decomposed into a rooted quadrangulation
counted by $y\frac{\pd}{\pd y}\Qb_j(x,y)$ and a sequence of
rooted quadrangulations counted by
$y\frac{\pd}{\pd y}\Qb_0(x,y)$.
Thus the generating function for $\net Q_1$ is
$$
\O\left(x^{-1}\frac{\pd \Qb_j(x,y)}{\pd y}
\Bigl(1-y\,\pd \Qb_0(x,y)/\pd y\Bigr)^{-1}\right).
$$

For convergence of $\sum(y\,\pd \Qb_0(x,y)/\pd y))^k$
it suffices to show that $y\,\pd \Qb_0(x,y)/\pd y<1$ for positive $x$
and $y$ since it is a power series with nonnegative coefficients.
Since
$$
\frac{y\,\pd \Qb_0(x,y)}{\pd y} ~=~ \frac{2(r+s)}{(2+r)(2+s)},
$$
the result is immediate.
Also note that this implies that $1-y\,\pd \Qb_0(x,y)/\pd y$ does
not vanish for $|x|\le\rho(r)$.

Similarly the quadrangulation $\map Q_2$ may contain contractible
2-cycles containing its root edge in its interior. So $\map Q_2$ is
decomposed into a rooted quadrangulation counted by $\Qb_{g-j}(x,y)$
and a sequence of rooted quadrangulations counted by
$y\frac{\pd}{\pd y}\Qb_0(x,y)$. Hence the generating function of the
quadrangulations with a separating non-contractible 2-cycle is
bounded above coefficient-wise by
\be \label{QbNC2}
\sum_{j=1}^{g-1}x^{-1}\Bigl(1-y\,\pd \Qb_0(x,y)/\pd y\Bigr)^{-2}\; \frac{\pd
\Qb_j(x,y)}{\pd y}\Qb_{g-j}(x,y),
\ee
which is algebraic with nonnegative coefficients.

Since $1-y\,\pd \Qb_0(x,y)/\pd y\ne 0$, the function given
in~(\ref{QbNC2}) has only one singularity on the circle of
convergence and near that singularity is $O((1-x/\rho(r))^p)$ where
$$
p = \left(\frac{3-5j}{2}-1\right) + \frac{3-5(g-j)}{2} ~=~
\frac{3-5g}{2}+\frac{1}{2}.
$$
Thus we can apply Lemma~\ref{Lemma:LLT} to see that the error is
negligible.

Next we consider quadrangulations counted by $\Qb_g(x,y)$ which
contain a non-separating non-contractible cycle $C$ of length $2$.
Cutting through $C$ gives a near-quadrangulation of genus $g-1$ with
two 2-cycles. After closing the resulting two 2-cycles, we obtain a
rooted quadrangulation $\map Q$ with two distinguished edges. The
quadrangulation $\map Q$ may contain contractible 2-cycles which
contain a distinguished edge in its interior.
Hence $\map Q$ is decomposed into a rooted quadrangulation counted by
$y^2\frac{\pd^2\Qb_{g-1}(x,y)}{(\pd y)^2}$ and two sequences of
rooted quadrangulations counted by $y\frac{\pd\Qb_0(x,y)}{\pd y}$.
Hence the bound in this case is
$$
\left(1-y\pd \Qb_0(x,y)/\pd y\right)^{-2} y^2
\frac{\pd^2\Qh_{g-1}(x,y) }{(\pd y)^2}.
$$
Reasoning as in the previous paragraph,
this gives a negligible contribution to the asymptotics.

Now Theorem~1(ii) follows from Lemma~1 and Theorem~3 using
$$
\frac{m}{n} = \frac{-d\,\log\rho(r)}{d\,\log
\eta_2(r)} = \frac{1+r}{r}
~~~\hbox{and}~~~
\si_2^2(r) = \frac{-d^2\,\log\rho(r)}{(d\,\log
\eta_2(r))^2} = \frac{(2+r)(1+2r)}{6r^2(1+r)}.
$$

\noindent{\bf Proof of Theorem~\ref{Thm:maps}(iii):} We prove that
almost all quadrangulations counted by $\QS(x,y)$ have no
non-contractible cycles of length 2 or 4. The argument is similar to
the one used above, and is identical to the one used in
\cite{BGRW96}. We note here
$$
\frac{m}{n} = \frac{-d\,\log\rho(r)}{d\,\log\eta_3(r)}
= \frac{3(1+r)}{1+2r}
~~~\hbox{and}~~~
\si_3^2(r) = \frac{-d^2\,\log\rho(r)}{(d\,\log \eta_3(r))^2}
= \frac{3r(2+r)}{2(1+r)(1+2r)^2}.
$$

\section{Face Widths of 3-Connected Maps and Graphs}\label{Sec:FW}

Robertson and Vitray~\cite{RV90} have shown that, if a
3-connected map $\map M$ in a surface $\Sigma_g$ of genus $g$
has $\fw(M)>2g+2$, then its underlying graph has
a unique embedding in $\Sigma_g$ and is not embeddable
in a surface of lower genus.

Our goal is to prove Theorem~\ref{Thm:fw} below. Then
Theorem~\ref{Thm:graphs}(i) follows from Theorem~\ref{Thm:maps} by
counting vertex-labeled, 3-connected, rooted maps. To obtain
Theorem~\ref{Thm:graphs}(iii) for 3-connected graphs, it suffices to
use (\ref{Qasymp}) with $r$ chosen so that $y=1$; that is,
$\eta_3(r)=1$. In other words, $r=\sqrt{7}/2-1$. This gives
$$
x_3=\r(r)=\frac{7\sqrt{7}-17}{32}\doteq 0.04751.
$$
By Comment~4 after Lemma~\ref{Lemma:LLT}, the number of edges is
concentrated around its mean which is asymptotically $\frac{3(1+r)}{1+2r}n$.

Applying the ``transfer theorem''~\cite[Sec.\thinspace VI.3]{FSbook}
to (\ref{Qasymp}) and using Theorem~\ref{Thm:fw}, one obtains
$$
\left(\frac{3(1+r)}{1+2r}n\right)\frac{G_g(n;3)}{n!} ~\sim~
\frac{C(r)n^{5(g-1)/2}}{4\,\Gamma\!\left(\frac{5g-3}{2}\right)}.
$$
After some algebra we obtain Theorem~\ref{Thm:graphs}(iii) for
3-connected graphs, with
\beqn
\b_3&=&\frac{2\sqrt{3}\,(1+2r)^2(1+r)^{3/2}(2+r)^{5/2}}{r^6}
~\doteq~1.48590\cdot 10^{5},\\
\a_3 &=& \frac{1}{4\b_3} ~\doteq~ 1.68248\cdot 10^{-6}. \eeqn

\begin{thm}[Large Face Width]\label{Thm:fw}
Fix $g>0$.
Let $L_g(x,y)=\sum_{n,m}L_g(n,m;c)x^ny^m$ where $L_g(n,m;c)$ is
the number of maps counted by $\Mh_g(n,m;3)$ that have face width
at least $c$.
Then
\beq
L_{g}(x,y)&=&x\QS_g(x,y)+\O(B_1(x,y)), \label{Lg}\\
\frac{\pd G_{g,3}(x,y)}{\pd y}
&=&\frac{x}{4y}\QS_g(x,y)+\O(B_2(x,y)),\label{GM3}
\eeq
where every singularity of $B_i$ is a singularity of $\QS_g$ and
$$
B_i(x,y) ~=~ O\left((1-x/\r(r))^{5(g-1)/2+1/2} \right)
\hbox{ ~as~ $x\to \r(r)$}
$$
for $y=\eta_3(r)$, uniformly for $r\in N(\epsilon)$.
\end{thm}

We show that almost all simple quadrangulations have no
non-contractible cycles of length less than any constant $c$.
We need only consider cycles of length $2k$ where $c\ge 2k>4$
since we may limit attention to simple quadrangulations.
Let $C$ be a non-contractible cycle of length $2k$ in a simple
quadrangulation counted by $\QS_g(x,y)$.
As in previous arguments, we consider separating and non-separating
separately

\medskip\noindent {\bf Case 1}.
Suppose $C$ is separating.
Cutting through $C$ and filling the two holes with
disks, we obtain a rooted simple quadrangulation $\map Q_1$ with a
distinguished face of degree $2k$, which has genus $0<j<g$, and
another rooted simple near-quadrangulation $\map Q_2$ with genus $g-j$
and root face degree $2k$.
We may quadrangulate the faces of degree $2k$ by inserting a
vertex in the interior of the face, but this may create separating
quadrangles near the cycle $C$.
We can get around this technical problem by gluing a special
near-quadrangulation $\map M_0$ to the face bounded by $C$.
For example, the near-quadrangulation $\map M_0$ can be constructed
using two copies of the $2k$-cycle, one inside the other,
adding edges between the two corresponding vertices of the cycles,
and inserting a new vertex inside the interior $2k$-cycle and
joining this new vertex to every other vertex of the cycle.
As a result we obtain a simple quadrangulation of genus $j$ with a
distinguished $\map M_0$, and another simple quadrangulation of genus
$g-j$ rooted at $\map M_0$.
Thus the generating function of simple quadrangulations in this case
is bounded by
$$
\O\left(x^{i}y^l\left(\sum_{j=1}^{g-1}\QS_{g-j}(x,y)
\frac{\pd \QS_{j}(x,y)}{\pd x}\right)\right)
$$
for some fixed integers $i,l$. As in previous arguments, this leads
to a negligible contribution.

\medskip\noindent {\bf Case 2}.
Now suppose $C$ is non-separating.
Cutting through $C$, filling the two holes with
disks, and then quadrangulating the resulting two faces as in
Case~1, we obtain a rooted simple quadrangulation of genus $g-1$
with two distinguished $\map M_0$.
Thus the generating function of simple quadrangulations in this
case is bounded by
$$\O\left(x^iy^l\frac{\pd^2 \QS_{g-1}(x,y)}{(\pd x)^2}\right)
$$
for some fixed integers $i,l$.
Again, the contribution is negligible.
This gives (\ref{Lg}).
Robertson and Vitray's result~\cite{RV90} implies that
$L_g(n,m;2g+3)n!/(4m)$ counts 3-connected graphs of genus $g$ with
face width at least $2g+3$ and so (\ref{GM3}) follows.

\section{From 3-connected graphs to 2-connected graphs}\label{Sec:3to2}

Since the results for 2-connected planar graphs follow
from~\cite{BGW02}, we assume $g>0$ in this section.

\begin{defn}[(Planar networks]
A {\em planar network} is a graph $\graph G$ together with two
distinguished vertices $v_0$ and $v_1$ (the {\em poles}) such that
the graph obtained by adding the edge $e=\{v_0,v_1\}$ (if it is not
already in $\graph G$) is 2-connected and planar. In contrast to the
usual labeled graph, the poles of a labeled network are not labeled.
\end{defn}

As in \cite{BGW02} we use $D(x,y)$ to denote the generating function
for planar networks. Thus $[(x^i/i!)y^m]\,D(x,y)$ is the number of
planar networks with $m$ edges $i$~vertices {\em not including} the
poles $v_0$ and $v_1$.

We will be expanding various functions about singularities. To help
us remember which coefficient goes with which function, we introduce
some notation. If $F(x)$ has a singularity at $x=r$ and we expand it
in powers of $(1-x/r)$, then $\cf F t$ denotes the coefficient of
$(1-x/r)^t$ in the expansion.

We begin with a review of some results for planar graphs. It is
convenient to use essentially the same notation and parametrization
as in \cite{BGW02}. That paper has three parameters, $u$, $v$ and
$t$. The parameters $u$ and $v$ are related to $r$ and $s$ by \be
u=\frac{r(2+s)}{4-rs}~~\mbox{and}~~v=\frac{s(2+r)}{4-rs};
\label{uvrs} \ee or equivalently, \be
r=\frac{2u}{1+v}~~\mbox{and}~~s=\frac{2v}{1+u}. \label{rsuv} \ee The
parameter $t$ is used on the singular curve $rs=1$ and is given by
$$
t=\frac{1}{1+2r}.
$$
It also uses the following functions of $t$. (When our notation
differs from~\cite{BGW02}, we have indicated the~\cite{BGW02}
notation parenthetically.)
\def\WasCalled#1#2{
  \hbox to3in{$\displaystyle #2$\hfill}
  \mbox{(called $#1$ in \cite{BGW02})}
}
\beqn
\a(t)&=&144 + 592t + 664t^2 + 135t^3 + 6t^4 - 5t^5 \\
\vphantom{\bigm|}
\b(t)&=&3t(1+t)(400 + 1808t + 2527t^2 + 1155t^3 + 237t^4 + 17t^5)\\
\gamma(t)&=&1296 +
10272\,t + 30920\,t^{2} + 42526\,t^{3} + 23135\,t^{4} \\
& &\BOp - 1482\,t^{5} - 4650\,t^{6} - 1358\,t^{7} - 405\,t^8 - 30t^9\\
\vphantom{\Biggm|_{\bigm|}^{\bigm|}}
h(t)&=& \frac {t^2(1-t)\,(18 + 36\,t + 5t^2)}{2(3+t)(1+2t)
\,(1+3\,t)^2}\\
\r_2(t)&=&\WasCalled{x_0}{\frac{(1+3t)(1-t)^3}{16t^3}}\\
\vphantom{\biggm|^{\bigm|}}
\l_2(t)&=&\WasCalled{y_0}{\frac {1+2t}{(1+3t)\,(1-t)}e^{-h(t)}-1}\\
\vphantom{\Biggm|_{\bigm|}^{\bigm|}}
\mu(t)& = &{(1+t) (3+t)^2 (1+2t)^2 (1+3t)^2 \l_2(t)\over
t^3(1+\l_2(t))\a(t)}\\
\si^2(t)&=&{(3+t)^2(1+2t)^2(1+3t)^2\l_2(t)\over
3t^6(1+t)(1+\l_2(t))^2\alpha(t)^3}\nonumber\\
& &\BOp\times
  \Bigl(3t^3(1+t)^2\a(t)^2-(1-t)(3+t)(1+2t)(1+3t)^2\l_2(t)\gamma(t) \Bigr)\\
\vphantom{\Biggm|_{\bigm|}^{\bigm|}}
\cf{\Do}0(t)&=& \WasCalled{D_0}{\frac{3t^2}{ (1-t)(1+3t)}}\\
\cf{\Do}1(t)&=&\WasCalled{D_2}{-\frac{48t^2(1+t)(1+2t)^2(18+6t+t^2)}{ (1+3t)\b(t)}}\\
\vphantom{\bigm|^{\bigm|}} \cf{\Do}{3/2}(t)&=&\WasCalled{D_3}{384
t^3 (1+t)^2 (1+2t)^2 (3+t)^2 \a(t)^{3/2}\beta(t)^{-5/2}}. \eeqn As
was pointed out in \cite{GN09}, a factor of $t$ is missing in $D_2$
of \cite{BGW02}. We note that $\r_2(t)=\r(r)$.

Throughout the rest of the paper, we adopt the following notation,
with $\e>0$ not necessarily the same at each appearance,
$$
T(\e) ~=~ \{te^{i\theta}:\e \!\le\! t \!\le\! 1-\e,\; |\theta|\le\e\}
~~~\hbox{and}~~~ \Delta(\r,\e)~=~\Domain.
$$
It is known \cite{BGW02,GN09} that for each $t\in T(\e)$,
$\Do(x,\l_2(t))$ and $G_{0,2}(x,\l_2(t))$ are all analytic in a
$\Delta(\r_2(t),\e)$ region. Also from~\cite{BGW02,GN09}, we have
\beq \Do(x,y)
&=&\cf{\Do}0(t)+\cf{\Do}1(t)(1-x/\r_2(t))+\cf{\Do}{3/2}(t)(1-x/\r_2(t))^{3/2}\label{D0}\\
& &\BOp+O\left((1-x/\r_2(t))^2\right),\nonumber \\
\vphantom{\Biggm|^{\Bigm|}_{\Bigm|}} \frac{\pd \Do}{\pd y}
&=&\frac{\cf{\Do}0'(t)}{\l_2'(t)}+\frac{\cf{\Do}1(t)\r_2'(t)}{\r_2(t)\l_2'(t)}
+O\left((1-x/\r_2(t))^{1/2}\right),\label{pD0}
 \eeq
as $x\to \r_2(t)$, uniformly for $y=\l_2(t)$ and $t\in T(\e)$.

We now turn our attention to $G_{g,2}(x,y)$ and $G_g(n,m;2)$.
Since the planar case $g=0$ has already been done \cite{BGW02,GN09}, we deal with
the nonplanar case and prove the following theorem.
\begin{thm} \label{G2asy}
Let $B_g(t)$ be as in Theorem~2(ii). There are generating functions
$E_{g,2}(x,y)$ which are analytic in a $\Delta(\r_2(t),\e)$ region
for each $t\in T(\e)$ such that \beqn G_{1,2}(x,y)&=&B_1(t)
\ln\left(\frac{1}{1-x/\r_2(t)}\right)+\O(E_{1,2}(x,y))\\
G_{g,2}(x,y)&=&B_g(t)
\,\Gamma\!\left(\frac{5g-5}{2}\right)\Bigl(1-x/\r_2(t)\Bigr)^{-5(g-1)/2}
\!\! +\O(E_{g,2}(x,y))~~~\mbox{for $g>1$.}
\eeqn
The radius of convergence $R(c)$ of $E_{g,2}(x,c)$ satisfies $R(c)>R(|c|)$ for
$c\ne |c|$. As $x\to \r_2(t)$, we have, uniformly for $y=\l_2(t)$
and $t\in T(\e)$,
$$E_{g,2}(x,y)= h(y)+O\left((1-x/\r_2(t))^{-5g/2+3}\right)$$
for some function $h(y)$.
\end{thm}

 \proof Since the planar case has been
done in~\cite{BGW02}, we will use induction on $g$ and assume $g>0$
below. Write $G_{g,2}(x,y) = F(x,y) + E(x,y)$ where $F(x,y)$ counts
2-connected graphs containing a unique nonplanar 3-connected
component and $E(x,y)$ counts the remaining 2-connected graphs. We
will analyze $F(x,y)$ and show that the contribution of $E(x,y)$ is
negligible.

The dominant singularity is extracted from the $F(x,y)$ part and the
remainder, along with the $E(x,y)$ bound, can be incorporated into
$E_{g,2}$.

We begin with $F$. A 2-connected graph $\net F$ counted by $F$
contains a unique 3-connected component of genus $g$ and all other
3-connected components of $\net F$ are planar.

Thus we have $$ F(x,y) ~=~ G_{g,3}(x,\Do(x,y)). $$ It follows from
(\ref{GM3}) that $$\frac{\pd }{\pd y}F(x,y) ~=~
\frac{x\QS_{g}(x,\Do(x,y))}{4\Do(x,y)}\frac{\pd \Do(x,y)}{\pd
y}+\O\left(B_2(x,\Do(x,y))\frac{\pd \Do(x,y)}{\pd y} \right),
$$ and hence
 \be \label{F}
F(x,y)=\int \frac{x\QS_{g}(x,\Do(x,y))}{4\Do(x,y)}\frac{\pd
\Do(x,y)}{\pd y}dy+\O \left(\int B_2(x,\Do(x,y))\frac{\pd
\Do(x,y)}{\pd y} dy\right).
 \ee
 Although we do not know $x\QS_{g}(x,y)$ exactly, we can still
obtain an asymptotic estimate for the above integral because the
coefficients of $\Do$ are nonnegative and we have (\ref{Qasymp}). We
first use Theorem~\ref{Thm:fw} and (\ref{Qasymp}) to obtain the
singular expansion for $x\QS_g(x,\Do(x,y))$ at the singularity
$x=\r(r)=\r_2(t)$, with $y=\l_2(t)$ fixed.  We have from
(\ref{Qasymp})
$$x\QS_g(x,\Do)= C(r)(1-x/\r(r))^{(3-5g)/2}+O\left((1-x/\r(r))^{(4-5g)/2}  \right),$$
as $x\to \r(r)$. As in the proofs of (\ref{Qasymp}) for $\Qh_g(x,y)$
and $\QS_g(x,y)$, it is important to note that $\Do=\Do(x,y)$ is a
function of $x$ for each fixed $y$, and hence $\r(r)$ is a function
of $x$ through the relation $\Do=\eta_3(r)$. It follows from
(\ref{D0}) that
 \beqn \left.\frac{d }{dx}
\left(1-\frac{x}{\r(r)}\right)\right|_{x=\r(r)}&=&
\frac{-1}{\r(r)}\left(1-\frac{\r'(r)}{\eta_3'(r)}\left.\frac{\pd \Do}{\pd x}\right|_{x=\r(r)}\right)\\
&=&\frac{-1}{\r(r)}\left(1+\frac{\r'(r)}{\eta_3'(r)}\frac{\Do^{[1]}}{\r(r)}\right)\\
&=& \frac{-1}{\r_2(t)}\frac{3(1+t)(1+3t)\a(t)}{\b(t)}. \eeqn
Hence

\beqn (1-x/\r(r))^{(3-5g)/2}&=&
\left(\frac{3(1+t)(1+3t)\a(t)}{\b(t)}\right)^{(3-5g)/2}(1-x/\r_2(t))^{(3-5g)/2}\\
& &+O\left((1-x/\r(r))^{(4-5g)/2} \right), \eeqn
as $x\to \r_2(t)$ with $y=\l_2(t)$ fixed. We remind the reader that $t$ and $r$ are
related by $t=\frac{1}{1+2r}$. Thus, temporarily using the notation
$$
H(t)=\sqrt{\frac{3}{1+r}}\frac{A_g(r)}{4\Do(x,y)(1+2r)}(2+r)^{5g-4}\left(\frac{3(1+t)(1+3t)\a(t)}{\b(t)}\right)^{(3-5g)/2}
\Gamma\left(\frac{5g-3}{2}\right),$$
we have
\beqn
F(x,y) &=& \int H(t)\Bigl(1-x/\r_2(t)\Bigr)^{(3-5g)/2}\,\frac{\pd \Do(x,y)}{\pd y}\frac{\l'_2(t)}{\r'_2(t)}d\r_2\\
& &\BOp+\O \left(\int B_2(x,\Do(x,y))\frac{\pd \Do(x,y)}{\pd
y}\frac{\l'_2(t)}{\r'_2(t)} d\r_2\right). \eeqn Noting that
$B_2(x,\Do(x,y))$ has a singular expansion at $x=\r_2(t)$ of lower
order, we obtain from (\ref{pD0}) that \beq
 F(x,y) &=&H(t)\left(\frac{\cf{\Do}0'(t)}{\l'_2(t)}+
 \frac{\r'_2(t)\cf{\Do}1(t)}{\l'_2(t)\r_2(t)}\right)\frac{\l'_2(t)\r_2(t)}{\r'_2(t)}f(\r_2)\nonumber\\
 & &\BOp+O \left(\Bigl(1-x/\r_2\Bigr)^{(6-5g)/2}\right)\nonumber\\
 &=&B_g(t)f(\r_2)+O\left( \Bigl(1-x/\r_2\Bigr)^{(6-5g)/2}\right),
 \eeq
where $B_g(t)$ is defined in Theorem~2(ii), $f_1(\r_2) =
-\ln(1-x/\rho_2)$, and
$$
f_g(\r_2) ~=~ \frac{(1-x/\r_2(t))^{-5(g-1)/2}}{5(g-1)/2}
~~\hbox{when $g>1$.}
$$

We now show that $E(x,y)$ is negligible compared with $F(x,y)$.

For each graph counted by $E(x,y)$, there are at least two nonplanar
3-connected components. In this case there is a 2-cut $\{a,b\}$ that
either splits $\graph G$ into two nonplanar pieces or gives a single
piece with a lower genus. We consider these two cases separately. As
in Section~\ref{Sec:Maps}, there is a non-contractible simple closed
curve $C$ intersecting $\graph G$ only at $a$ and $b$. As an aside,
we note that this means the face width of $\graph G$ is at most 2
and hence intuitively the graphs in this class should be negligible;
however, we have not proved a large face-width result for
2-connected graphs. The following analysis basically proves such a
large face-width result and is very similar to the one used above
for 3-connected graphs (maps).

\medskip\noindent {\bf Case 1}.
Cutting through $C$ splits $\graph G$ into 2-connected graphs
$\graph G_1$ and $\graph G_2$ such that $\graph G_1$ is embeddable
in the orientable surface of genus $j>0$ and $\graph G_2$ is
embeddable in the orientable surface of genus $g-j>0$. Also $\graph
G_1$ and $\graph G_2$ each has a distinguished edge (joining
vertices $a$ and $b$). Hence the generating function of the
2-connected graphs in this case is bounded by (applying
Lemma~\ref{Lemma:O(GF)})
$$
\sum_{j=1}^{g-1}\O\left(\frac{\pd G_{g-j,2}(x,y)}{\pd y}\frac{\pd
G_{j,2}(x,y)}{\pd y}\right).
$$

\smallskip\noindent {\bf Case 2}.
Cutting through $C$ reduces $\graph G$ into a 2-connected graph
$\net G_1$ which is embeddable in the orientable surface of genus
$g-1$, and $\net G_1$ has two distinguished edges (joining the
copies of $a$ and $b$). Hence the generating function in this case
is bounded by
$$
\O\left(\frac{\pd^2 G_{g-1,2}(x,y)}{(\pd y)^2}\right).
$$

By induction, it is easily seen that the contributions in both cases
satisfy Lemma~\ref{Lemma:LLT} with the same parameters as $F(x,y)$,
except that the exponent of $n$ obtained in the asymptotics is less
than the exponent of $n$ in the asymptotics for $F$.

We need to establish Lemma~\ref{Lemma:LLT}(a,e). It is important to
note that the dominant singularities of $G_{g,3}(x,\Do(x,y))$ are
the same for each genus $g$ because $G_{g,3}(x,y)$ (more precisely
$\QS_g(x,y)$) have the same dominant singularities.

This completes the proof of Theorem~\ref{G2asy}.

\medskip

Now Theorem~\ref{Thm:graphs}(ii) follows immediately using
Lemma~\ref{Lemma:LLT}.

Theorem~\ref{Thm:graphs}(iii) for 2-connected graphs follows by
setting $y=\l_2(t)=1$ in Theorem~5 (i.e., $t=\tO\doteq 0.62637$) and
applying the ``transfer'' theorem.  We note that \beqn \b_2& =&
\frac{8}{9(1+\tO)(1-\tO)^6}\left(\frac{\b(\tO)}{\a(\tO)}\right)^{5/2}
~\doteq~7.6150\cdot 10^4 ,\\
\a_2 &=& \frac{1}{4\b_2}~\doteq~ 3.28299\cdot 10^{-6}.
 \eeqn

\section{From 2-connected graphs to 1-connected graphs}\label{Sec:2to1}

Since the composition depends only on the vertices, there is no need
to keep track of the number of edges if we only care about the
number of graphs with $n$ vertices. This makes the arguments much
simpler as we are dealing with univariate functions. From now on, we
will focus on $y=1$, although the results extend to all $y$ near $1$
as done in \cite{GN09} for planar graphs. We note that would be
possible to extend the result to the whole range of $y$, provided
that the condition $R(y)>R(|y|)$ when $y\ne |y|$ for the radius of
convergence in Lemma~\ref{Lemma:LLT} can be verified for
$G_{g,1}(x,y)$. However, we have not verified this technical
condition.

Since the planar case is dealt with in \cite{GN09}, we assume $g>0$.

Let $x_1$ be the smallest positive singularity of $G_{0,1}(x)$.
Gim\'enez and Noy~\cite[p.\thinspace320]{GN09} showed that
\be\label{x12} x_2 = x_1G_{0,1}'(x_1) \ee and $G_{0,1}(x)$ is
analytic in a $\Delta(x_1,\e)$ region.

As in the previous section, let $\tO\doteq 0.62637$ be determined by
$\l_2(\tO)=1$. From \cite[Lemma 6]{GN09}, we have the following
singular expansion at $x_2=\r_2(\tO)\doteq 0.03819$, \be
\label{G2planar} G_{0,2}(x) ~=~ \cf{G_{0,2}}0 +
\cf{G_{0,2}}1(1-x/x_2) + \cf{G_{0,2}}2(1-x/x_2)^2 +
\cf{G_{0,2}}{5/2}(1-x/x_2)^{5/2} + \ldots, \ee where
$\cf{G_{0,2}}j=\cf{G_{0,2}}j(\tO)$, and in particular
$$
\cf{G_{0,2}}0\doteq 7.397\cdot 10^{-4},~~
\cf{G_{0,2}}1\doteq -1.4914\cdot 10^{-3}
~~\hbox{and}~~
\cf{G_{0,2}}2\doteq 7.672\cdot 10^{-4}.
$$
Define
\beqn
A&=&
\frac{(3\tO-1)(1+\tO)^3\ln(1+\tO)}{16\tO^3}+\frac{(1+3\tO)(1-\tO)^3\ln(1+2\tO)}{32\tO^3}\\
&&\BOp+\frac{(1-\tO)(185\tO^4+698\tO^3-217\tO^2-160\tO+6)}{64\tO(1+3\tO)^2(3+\tO)},\\
x_1&=& {\textstyle\frac{1}{16}}\sqrt{1+3\tO} (1-\tO)^3\tO^{-3}e^A
~\doteq~ 0.03673. \eeqn It was shown in \cite{GN09} that \beq
\label{G1planar} G_{0,1}(x) &=& \cf{G_{0,1}}0+\cf{G_{0,1}}1(1-x/x_1)
+\cf{G_{0,1}}2(1-x/x_1)^2+\cf{G_{0,1}}{5/2}(1-x/x_1)^{5/2}+\ldots\qquad\\
P(x) &:=& xG'_{0,1}(x) ~=~ \cf P0+\cf P1(1-x/x_1)+\cf
P{3/2}(1-x/x_1)^{3/2}+\ldots \label{Px},
 \eeq
where
$$
\cf P0=-\cf{G_{0,1}}1, ~~\cf P1=-2\cf{G_{0,1}}2-\cf{G_{0,1}}0\doteq
-0.03979
~~\hbox{and}~~
\cf P{3/2}=-5\cf{G_{0,1}}{5/2}/2.
$$
We also note that \cite[(4.7)]{GN09}
$\cf{G_{0,1}}0=G_{0,1}(x_1)=x_2+\cf{G_{0,2}}0+\cf{G_{0,2}}1\doteq
0.03744$. The following theorem summarizes the main results of this
section.
\begin{thm} \label{Thm:G1asy}
Fix $g>0$. We have $G_{g,1}(x)= F(x) + \O(E(x))$ where
\begin{itemize}
\item[(i)]
$F(x)$ and $E_g(x)$ are analytic in a $\Delta(x_1,\e)$;
\item[(ii)]
as $x\to x_1$,
$$
F(x) ~\sim~ \cases{ \displaystyle \a_2\b_2t_1
\ln\left(\frac{1}{1-x/x_1}\right)_{\vphantom{|}} &if $g=1$,\cr
\displaystyle
\a_2\b_2^gt_g\,\Gamma\!\left(\frac{5g-5}{2}\right)^{\vphantom{|}}
\left(\frac{-x_2}{\cf P1}\right)^{5(g-1)/2}(1-x/x_1)^{-5(g-1)/2} &if
$g>1$;\cr}
$$
\item[(iii)]
as $x\to x_1$, $E_1(x) = C+O\left((1-x/x_1)^{1/2}\right)$ for some
constant $C$ and \hfill\break $E(x) =
O\left((1-x/x_1)^{-5g/2+3}\right)$ when $g>1$.
\end{itemize}
\end{thm}

\proof We again apply induction on $g$.
Let $\graph G$ be
a connected graph of genus $g$ rooted at a vertex $v$.
It is well known that $\graph G$ is (uniquely) decomposed into
a set of blocks (2-connected pieces) and the genus of $\graph G$
is the sum of the genera of all blocks~\cite{BHKY62}.
We divide all connected graphs of genus $g>0$ into two
classes according to whether there is a block of genus $g$ or not
and will show that the second class is negligible.

\medskip\noindent{\bf Case 1 (Genus $g$ block)}.
We attach a planar 1-connected graph to each vertex of the genus $g$ to
connected block.
Thus the generating function for this case is
$$
F(x) ~=~ G_{g,2}(xG'_{0,1}(x)).
$$
Since $G_{g,2}(x)$ is bounded termwise above and below by functions
analytic in a $\Delta(x_2,\e)$ region, it follows from (\ref{x12}),
the same holds for $F(x)$ in a $\Delta(x_1,\e)$ region. For $g>1$,
it follows from Theorem~5 and (\ref{Px}) that \beqn G_{g,2}(x) &=&
\a_2\b_2^gt_g\,\Gamma\!\left(\frac{5g-5}{2}\right)
(1-x/x_2)^{-5g/2+5/2}+O\left((1-x/x_2)^{-5g/2+3}\right)\\
G_{g,2}(xG'_{0,1}(x))&=&
\a_2\b_2^gt_g\,\Gamma\!\left(\frac{5g-5}{2}\right)
(-\cf P1/x_2)^{-5g/2+5/2}(1-x/x_1)^{-5g/2+5/2}\\
& &\BOp+ O\left((1-x/x_1)^{-5g/2+3}\right). \eeqn It follows that
\beq\nonumber F(x)&=&
\a_2\b_2^gt_g\,\Gamma\!\left(\frac{5g-5}{2}\right)
(-\cf P1/x_2)^{-5g/2+5/2}(1-x/x_1)^{-5g/2+5/2}\\
& &\BOp+ O\left((1-x/x_1)^{-5g/2+3}\right). \label{F(x) asymp}
\eeq
The formula for $g=1$ is similar except that it involves a logarithm:
\be\label{F1(x) asymp}
F(x) ~=~
\a_2\b_2 t_1 \ln\left(\frac{1}{1-x/x_1}\right) + O(1).
\ee
As in previous proofs, we use bounds on the functions to split
$F$ into $F_g$ and a contribution to $E_g$ which are analytic in
a $\Delta(x_1,\e)$.

\smallskip\noindent{\bf Case 2 (No genus $g$ block)}.
In this case, there is at least one vertex $v$ such that $\graph G$
can be viewed as two nonplanar graphs joined at $v$.
Hence an upper bound for graphs in this class is given by
the generating function
$$
\sum_{j=1}^{g-1}xG'_{j,1}(x)G'_{g-j,1}(x).
$$
It follows by induction on $g$ that each summand is bounded
by a function analytic in a $\Delta(x_1,\e)$ region and,
as $x\to x_1$ in this region each bound is bounded by
\be\label{no g}
O\left((1-x/x_1)^{-5j/2+3/2}(1-x/x_1)^{-5(g-j)/2+3/2}\right) ~=~
O\left((1-x/x_1)^{-5g/2+3}\right).
\ee
This completes the proof of Theorem~\ref{Thm:G1asy}.

\medskip
Now Theorem~\ref{Thm:graphs}(iii) for 1-connected graphs follows
immediately using the ``transfer theorem''. We obtain
$$
\b_1 ~=~ \left(\frac{-x_2}{\cf P1}\right)^{\!5/2}\!\!\!\b_2 ~\doteq~
6.87242\cdot 10^4,~~~\hbox{and}~~~ \a_1 ~=~ \frac{1}{4\b_1} ~\doteq~
3.63773\cdot 10^{-6}.
$$

\section{From 1-connected graphs to all graphs}\label{Sec:1toAll}

The case $g=0$ is treated in \cite{GN09}
We treat $g>1$.
The case $g=1$ is similar except that $\ln(1-x/x_1)$ appears.
Let $F(x)$ denote the generating function of these graphs containing
a connected component of genus $g$. Then we have
\beqn
F(x) &= &G_{g,1}(x)\exp(G_{0,1}(x))\\
&=& \a_1\b_1^g\exp(G_{0,1}(x_1))t_g\, \Gamma\!\left(\frac{5g-5}{2}\right)(1-x/x_1)^{5(1-g)/2} \\
&&\BOp+O\left((1-x/x_1)^{-5g/2+3}\right).
 \eeqn
Again the case that there are two components with positive genus is
(by induction) bounded by
$$
O\left(\sum_{j=1}^{g-1}G_{j,1}(x)G_{g-j,0}(x)\right)
~=~ O\left((1-x/x_1)^{-5g/2+3}\right).
$$
Thus
\beqn
G_{g,0}(x)
&=& G_{g,1}(x)\exp(G_{0,1}(x))+O\left((1-x/x_1)^{-5g/2+3}\right)\\
&=&
\a_1\b_1^{g}\exp(G_{0,1}(x_1))t_g\,\Gamma\!\left(\frac{5g-5}{2}\right)
  (1-x/x_1)^{5(1-g)/2}\\
&&\BOp+O\left((1-x/x_1)^{-5g/2+3}\right). \eeqn
This completes the proof of Theorem~\ref{Thm:graphs}
(using the ``transfer'' theorem again) with
$$
\a_0 ~=~\a_1\exp(G_{0,1}(x_1)) ~\doteq~ 3.77651\cdot 10^{-6}
~~\hbox{and}~~ \b_0=\b_1.
$$

\section{A formula for $A_g(r)$}\label{Sec:Ag(r)}

In this section we obtain a formula for $A_g(r)$ using~\cite{BCR93}
and recently derived information~\cite{G1x} for $t_g(r)$.

Let $T_g(n,j)$ be the number of rooted maps of genus $g$ with
$i$ faces and $j$ vertices.

By duality, we may interchange the role of vertices and faces, and we do so.
By Euler's formula, $T_g(n,j)$ is also the number of rooted maps of
genus $g$ with $i$ vertices and $m=j+n+2g-2$ edges.

By Theorem~\ref{Thm:maps}, we have
\beq\nonumber
T_g(n,j)&=&[x^ny^m]\Mh_g(x,y)\\
\label{T_g}
&\sim& \vphantom{\biggl|}
\Bigl(C_1(r)A_g(r)n^{5g/2-3}\Bigr)\rho(r)^{-n}\eta_1(r)^{-m}\\
\nonumber
&\sim&
C_1(r)A_g(r)(n/j)^{5g/4-3/2}(nj)^{5g/4-3/2}\rho(r)^{-n}\eta_1(r)^{-n-j+2-2g}
\eeq
Note that
$$
\frac{j}{n} ~=~ \frac{m}{n}-1 +\frac{2-2g}{n}
~=~ \frac{1+2r}{r^2(2+r)} + \frac{2-2g}{n}.
$$

It follows that the value of $r$ in \cite[Theorem~2]{BCR93}
differs from our $r$ by $O(1/n)$.
Replacing one $r$ with the other inside the large parentheses of (\ref{T_g}) does
not change the asymptotics.
We must show that is also true for $f=\rho(r)^{-n}\eta_1(r)^{-m}$.
This can be done by expanding $\log f$ in a power series about $r$ and noting that
the linear term vanishes after we set $m/n$ to the value given in Theorem~\ref{Thm:maps}(i).
It follows that replacing $r$ by $r+O(1/n)$ changes $\log f$ by $nO(1/n^2)=o(1)$.
Hence we may freely use either value of $r$ in (\ref{T_g}).
Thus we obtain
$$
T_g(n,j) ~\sim~
\left(C_1(r)A_g(r)\eta_1(r)^{2-2g}\left(\frac{r^2(2+r)}{1+2r}\,nj\right)^{5g/4-3/2}\right)
\left(\rho(r)\eta_1(r)\right)^{-n}\eta_1(r)^{-j}.
$$
Comparing this with \cite[Theorem~2]{BCR93} we obtain
$$
t_g(r) ~=~ C_1(r)A_g(r)\eta_1(r)^{2-2g}\left(\frac{r^2(2+r)}{1+2r}\right)^{5g/4-3/2}
$$
and so by Theorem~\ref{Thm:maps}(i)

\beq\nonumber
A_g(r)
&=& \frac{\eta_1(r)^{2g-2}}{C_1(r)}
\left(\frac{1+2r}{r^2(2+r)}\right)^{5g/4-3/2}_{\vphantom{\Bigm|}}t_g(r)\\
\nonumber
&=& \frac{2^4r^3(2+r)^{1/2}(1+r+r^2)^{7/2}(4+7r+4r^2)^{1/2}}{(1+2r)^3}\\
\label{Ag and tg}
&&\BOp\times\left(\frac{(1+2r)^{13/2}}{2^8r^5(1+r+r^2)^8(2+r)^{5/2}}\right)^{g/2}
\!t_g(r).
\eeq
It was shown by the second author~\cite{G1x} that
\be \label{tgr}
t_g(r)=c(r)[d(r)]^g\, t_g,
\ee
where
\beqn
c(r)&=& \frac{r^3(1+2r)(2+r)}{32\sqrt{\pi}(4+7r+4r^2)^{1/2}(1+r+r^2)^{7/2}_{\vphantom{\bigm|}}},\\
d(r)&=& \frac{32\sqrt{3}(1+r+r^2)^4(1+r)^{3/2}}{r^{7/2}(2+r)^{5/4}(1+2r)^{5/4}}.
\eeqn
Combining (\ref{Ag and tg}) and (\ref{tgr}) gives
$$
A_g(r)
~=~ \frac{r^6(2+r)^{3/2}}{2\sqrt\pi(1+2r)^2}
\left(\frac{12(1+2r)^4(1+r)^3}{r^{12}(2+r)^5}\right)^{g/2}\!t_g.
$$

\section{Remarks on Nonorientable Surfaces}

We believe the study of maps on nonorientable surfaces and graphs
embeddable in nonorientable surfaces proceeds in a manner akin to
the orientable case presented here, as happened for maps
in~\cite{BCR93}. In particular, Theorems 1 and 2 hold for the
nonorientable surface with $2g$ crosscaps, with $t_g$ replaced by
the nonorientable map asymptotics constants $p_g$ in the expression
of $A_g(r)$ in Theorem~1. The major differences to be expected are
as follows.
\begin{itemize}
\item
The projective plane will require some special care because the
singular expansion of the generating function behaves like
$$M_{1/2}(x,y)=g(y)+h(y)(1-x/\r(y))^{1/4}+\cdots
$$
whose dominant term has a positive exponent. Consequently, the
exponent of the dominant term in the singular expansion of the
product $M_{1/2}(x,y)M_{g-1/2}(x,y)$ is not simply the sum
$(1/4)+(3/2)-5(g-1/2)/2$ as indicated on page 13. Rather, it should
be $(3/2)-5(g-1/2)/2$ when $g>1$ and $1/4$ when $g=1$. It can still
be checked that the exponent of the singular expansion of a
product like (\ref{QbNC2}) is higher than $(3-5g)/2$.
\item
Careful attention to the proofs in this paper show that
relative errors in the asymptotics are $O(n^{-1/2})$.
In the nonorientable case they should be $O(n^{-1/4})$.
\item
Unlike $t_g$, no simple recursion is known for the nonorientable map
asymptotics constant $p_g$. However, Garoufalidis and Mari\~{n}o
\cite{GM09} conjectured that
$$
p_g ~=~ \frac{v_{2g-1}}{2^{g-2}\Gamma\left(\frac{5g-3}{2}\right)},
$$
where $v_g$ satisfies
$$
v_g ~=~ \frac{1}{2\sqrt{3}}\left(-3a_{g/2}+\frac{5g-6}{2}v_{g-1}+\sum_{k=1}^{g-1}v_kv_{g-k}\right),$$

and $a_j$ is defined in (\ref{ag}), with the understanding that $a_j=0$
when $j$ is not an integer.

It was shown in~\cite{G1x} that $p_g(r)$ and $p_g$ satisfy the same
relation as $t_g(r)$ and $t_g$ expressed in (\ref{tgr}), namely
$p_g(r)=c(r)[d(r)]^gp_g$.
\end{itemize}

\end{document}